\documentclass{article}
\usepackage{amssymb,amsmath,theorem,euscript}

\input{epsf}

\newcounter{sec}

\def\sm{\smallskip}

\newcounter{punct}[sec]

\def\punct{\refstepcounter{punct}{\arabic{sec}.\arabic{punct}.  }} 

\def\COUNTERS{\addtocounter{sec}{1}
              \setcounter{punct}{0}
          \setcounter{equation}{0}
          \setcounter{theorem}{0}
                  }

\newtheorem{theorem}{Theorem}[sec]

\newtheorem{lemma}[theorem]{Lemma}

\begin{document}

 \def\ov{\overline}
\def\wt{\widetilde}
\def\wh{\widehat}
 \newcommand{\rk}{\mathop {\mathrm {rk}}\nolimits}
  \newcommand{\sgn}{\mathop {\mathrm {sgn}}\nolimits}
\newcommand{\Aut}{\mathop {\mathrm {Aut}}\nolimits}
\newcommand{\Out}{\mathop {\mathrm {Out}}\nolimits}
 \newcommand{\tr}{\mathop {\mathrm {tr}}\nolimits}
  \newcommand{\diag}{\mathop {\mathrm {diag}}\nolimits}
  \newcommand{\supp}{\mathop {\mathrm {supp}}\nolimits}
  \newcommand{\indef}{\mathop {\mathrm {indef}}\nolimits}
  \newcommand{\dom}{\mathop {\mathrm {dom}}\nolimits}
  \newcommand{\im}{\mathop {\mathrm {im}}\nolimits}
   \newcommand{\ind}{\mathop {\mathrm {ind}}\nolimits}
    \newcommand{\codim}{\mathop {\mathrm {codim}}\nolimits}
 
\renewcommand{\Re}{\mathop {\mathrm {Re}}\nolimits}
\renewcommand{\Im}{\mathop {\mathrm {Im}}\nolimits}

\def\Br{\mathrm {Br}}

\def\SL{\mathrm {SL}}
\def\Diag{\mathrm {Diag}}
\def\SU{\mathrm {SU}}
\def\GL{\mathrm {GL}}
\def\U{\mathrm U}
\def\OO{\mathrm O}
 \def\Sp{\mathrm {Sp}}
 \def\SO{\mathrm {SO}}
\def\SOS{\mathrm {SO}^*}
 \def\Diff{\mathrm{Diff}}
 \def\Vect{\mathfrak{Vect}}
\def\PGL{\mathrm {PGL}}
\def\PU{\mathrm {PU}}
\def\PSL{\mathrm {PSL}}
\def\Symp{\mathrm{Symp}}
\def\End{\mathrm{End}}
\def\Mor{\mathrm{Mor}}
\def\Aut{\mathrm{Aut}}
 \def\PB{\mathrm{PB}}
 \def\cA{\mathcal A}
\def\cB{\mathcal B}
\def\cC{\mathcal C}
\def\cD{\mathcal D}
\def\cE{\mathcal E}
\def\cF{\mathcal F}
\def\cG{\mathcal G}
\def\cH{\mathcal H}
\def\cJ{\mathcal J}
\def\cI{\mathcal I}
\def\cK{\mathcal K}
 \def\cL{\mathcal L}
\def\cM{\mathcal M}
\def\cN{\mathcal N}
 \def\cO{\mathcal O}
\def\cP{\mathcal P}
\def\cQ{\mathcal Q}
\def\cR{\mathcal R}
\def\cS{\mathcal S}
\def\cT{\mathcal T}
\def\cU{\mathcal U}
\def\cV{\mathcal V}
 \def\cW{\mathcal W}
\def\cX{\mathcal X}
 \def\cY{\mathcal Y}
 \def\cZ{\mathcal Z}
\def\0{{\ov 0}}
 \def\1{{\ov 1}}
 \def\frA{\mathfrak A}
 \def\frB{\mathfrak B}
\def\frC{\mathfrak C}
\def\frD{\mathfrak D}
\def\frE{\mathfrak E}
\def\frF{\mathfrak F}
\def\frG{\mathfrak G}
\def\frH{\mathfrak H}
\def\frI{\mathfrak I}
 \def\frJ{\mathfrak J}
 \def\frK{\mathfrak K}
 \def\frL{\mathfrak L}
\def\frM{\mathfrak M}
 \def\frN{\mathfrak N} \def\frO{\mathfrak O} \def\frP{\mathfrak P} \def\frQ{\mathfrak Q} \def\frR{\mathfrak R}
 \def\frS{\mathfrak S} \def\frT{\mathfrak T} \def\frU{\mathfrak U} \def\frV{\mathfrak V} \def\frW{\mathfrak W}
 \def\frX{\mathfrak X} \def\frY{\mathfrak Y} \def\frZ{\mathfrak Z} \def\fra{\mathfrak a} \def\frb{\mathfrak b}
 \def\frc{\mathfrak c} \def\frd{\mathfrak d} \def\fre{\mathfrak e} \def\frf{\mathfrak f} \def\frg{\mathfrak g}
 \def\frh{\mathfrak h} \def\fri{\mathfrak i} \def\frj{\mathfrak j} \def\frk{\mathfrak k} \def\frl{\mathfrak l}
 \def\frm{\mathfrak m} \def\frn{\mathfrak n} \def\fro{\mathfrak o} \def\frp{\mathfrak p} \def\frq{\mathfrak q}
 \def\frr{\mathfrak r} \def\frs{\mathfrak s} \def\frt{\mathfrak t} \def\fru{\mathfrak u} \def\frv{\mathfrak v}
 \def\frw{\mathfrak w} \def\frx{\mathfrak x} \def\fry{\mathfrak y} \def\frz{\mathfrak z} \def\frsp{\mathfrak{sp}}
 
 \def\frgl{\mathfrak{gl}}
 \def\bfa{\mathbf a} \def\bfb{\mathbf b} \def\bfc{\mathbf c} \def\bfd{\mathbf d} \def\bfe{\mathbf e} \def\bff{\mathbf f}
 \def\bfg{\mathbf g} \def\bfh{\mathbf h} \def\bfi{\mathbf i} \def\bfj{\mathbf j} \def\bfk{\mathbf k} \def\bfl{\mathbf l}
 \def\bfm{\mathbf m} \def\bfn{\mathbf n} \def\bfo{\mathbf o} \def\bfp{\mathbf p} \def\bfq{\mathbf q} \def\bfr{\mathbf r}
 \def\bfs{\mathbf s} \def\bft{\mathbf t} \def\bfu{\mathbf u} \def\bfv{\mathbf v} \def\bfw{\mathbf w} \def\bfx{\mathbf x}
 \def\bfy{\mathbf y} \def\bfz{\mathbf z} \def\bfA{\mathbf A} \def\bfB{\mathbf B} \def\bfC{\mathbf C} \def\bfD{\mathbf D}
 \def\bfE{\mathbf E} \def\bfF{\mathbf F} \def\bfG{\mathbf G} \def\bfH{\mathbf H} \def\bfI{\mathbf I} \def\bfJ{\mathbf J}
 \def\bfK{\mathbf K} \def\bfL{\mathbf L} \def\bfM{\mathbf M} \def\bfN{\mathbf N} \def\bfO{\mathbf O} \def\bfP{\mathbf P}
 \def\bfQ{\mathbf Q} \def\bfR{\mathbf R} \def\bfS{\mathbf S} \def\bfT{\mathbf T} \def\bfU{\mathbf U} \def\bfV{\mathbf V}
 \def\bfW{\mathbf W} \def\bfX{\mathbf X} \def\bfY{\mathbf Y} \def\bfZ{\mathbf Z} \def\bfw{\mathbf w}
 \def\R {{\mathbb R }} \def\C {{\mathbb C }} \def\Z{{\mathbb Z}} \def\H{{\mathbb H}} \def\K{{\mathbb K}}
 \def\N{{\mathbb N}} \def\Q{{\mathbb Q}} \def\A{{\mathbb A}} \def\T{\mathbb T} \def\P{\mathbb P} \def\G{\mathbb G}
 \def\bbA{\mathbb A} \def\bbB{\mathbb B} \def\bbD{\mathbb D} \def\bbE{\mathbb E} \def\bbF{\mathbb F} \def\bbG{\mathbb G}
 \def\bbI{\mathbb I} \def\bbJ{\mathbb J} \def\bbK{\mathbb K} \def\bbL{\mathbb L} \def\bbM{\mathbb M} \def\bbN{\mathbb N} \def\bbO{\mathbb O}
 \def\bbP{\mathbb P} \def\bbQ{\mathbb Q} \def\bbS{\mathbb S} \def\bbT{\mathbb T} \def\bbU{\mathbb U} \def\bbV{\mathbb V}
 \def\bbW{\mathbb W} \def\bbX{\mathbb X} \def\bbY{\mathbb Y} \def\kappa{\varkappa} \def\epsilon{\varepsilon}
 \def\phi{\varphi} \def\le{\leqslant} \def\ge{\geqslant}

\def\UU{\bbU}
\def\Mat{\mathrm{Mat}}
\def\tto{\rightrightarrows}

\def\Gr{\mathrm{Gr}}

\def\graph{{graph}}

\def\O{\mathrm{O}}

\def\la{\langle}
\def\ra{\rangle}

\def\B{\mathrm B}
\def\Int{\mathrm{Int}}
\def\LGr{\mathrm{LGr}}


\def\I{\mathbb I}
\def\M{\mathbb M}
\def\T{\mathbb T}

\def\Lat{\mathrm{Lat}}
\def\LLat{\mathrm{LLat}} 
\def\Mod{\mathrm{Mod}}
\def\LMod{\mathrm{LMod}}
\def\Naz{\mathrm{Naz}}
\def\naz{\mathrm{naz}}
\def\bNaz{\mathbf{Naz}}
\def\AMod{\mathrm{AMod}}
\def\ALat{\mathrm{ALat}}
\def\MAT{\mathrm{MAT}}
\def\Mar{\mathrm{Mar}}

\def\Ver{\mathrm{Vert}}
\def\Bd{\mathrm{Bd}}
\def\We{\mathrm{We}}
\def\Heis{\mathrm{Heis}}
\def\Pol{\mathrm{Pol}}
\def\Ams{\mathrm{Ams}}
\def\Herm{\mathrm{Herm}}

\def\F{\bbF}
\def\V{\bbV}

\def\Exp{\mathsf{Exp}}

 \newcommand{\gr}{\mathop {\mathrm {gr}}\nolimits}

\def\principal{\mathrm{principal}}
\def\discrete{\mathrm{discrete}}
\def\tempered{\mathrm{tempered}}
 
 \def\kos{/\!\!/}
 
\begin{center}
 \Large \bf
 
 The Fourier transform on the group $\GL_2(\R)$  and 
 the action of the overalgebra $\frgl_4$

 \large \sc
 
 \bigskip
 
 Yu.A.Neretin%
 \footnote{Supported by the grant FWF, Project P28421.}
 
\end{center}

{\small We define a kind of 'operational calculus' for $\GL_2(\R)$.  Namely,
	 the group $\GL_2(\R)$ can be regarded as an open dense chart in the Grassmannian
	of 2-dimensional subspaces in $\R^4$. Therefore the group $\GL_4(\R)$ acts
	in $L^2$ on $\GL_2(\R)$. We transfer the corresponding action of the Lie algebra
	$\frgl_4$ to the Plancherel decomposition of $\GL_2(\R)$, 
	the algebra acts by differential-difference
	operators with shifts in an imaginary direction. We also write similar formulas
	for the action of $\frgl_4\oplus \frgl_4$ in the Plancherel decomposition of
	$\GL_2(\C)$.}

\section{The statement of the paper}

\COUNTERS

{\bf \punct The group $\GL_(\R)$.%
	\label{ss:R}}  
Let $\GL_n(\R)$ be the group of invertible real matrices 
of order $n$. We denote elements of $\GL_2(\R)$ by
$$
X=\begin{pmatrix}
x_{11}&x_{12}\\
x_{21}&x_{22}
\end{pmatrix}
.
$$
 The Haar measure on $\GL_2(\R)$ is given by
$$
(\det X)^{-2}\,dX=
(\det X)^{-2} \,dx_{11} \,dx_{12}\,\,dx_{21}\,dx_{22}.
$$

Recall some basic facts on representations of
the group $\GL_2(\R)$, for systematic exposition of the representation
theory of $\SL_2(\R)$, see, e.g., \cite{GGV}, \cite{HT}.
Denote by $\R^\times$ the multiplicative group of $\R$.
Let $\mu\in\C$ and $\epsilon\in\Z/2\Z$. We define the function
$x^{\mu\kos  \epsilon}$ on $\R^\times$ by
$$
x^{\mu\kos  \epsilon}:=|x|^\mu \sgn(x)^\epsilon.
$$
Denote by $\Lambda$
the set of all
collections
$$(\mu_1,\epsilon_1;\mu_2,\epsilon_2),$$
 i.e.,
 $$\Lambda\simeq\C\times \Z_2\times \C\times \Z_2.$$ 
For each element of $\Lambda$ we define a representation  
$T_{\mu,\epsilon}$ of $\GL_2(\R)$ in the space of functions 
on $\R$ by 
\begin{multline}
T_{\mu_1,\epsilon_1;\mu_2,\epsilon_2}\begin{pmatrix}
x_{11}&x_{12}\\
x_{21}&x_{22}
\end{pmatrix} \phi(t)=\\=
\phi\Bigl(\frac{x_{12}+t x_{22}}{x_{11}+t x_{21}}\Bigr)
\cdot(x_{11}+t x_{21})^{-1+\mu_1-\mu_2\kos  \epsilon_1-\epsilon_2} \det\begin{pmatrix}
x_{11}&x_{12}\\
x_{21}&x_{22}
\end{pmatrix}^{1/2+\mu_2\kos  \epsilon_2}
.
\label{eq:T-mu}
\end{multline}

A space of functions can be specified in various ways, it is convenient
to consider   the space $C^\infty_{\mu_1-\mu_2,\epsilon_1-\epsilon_2}$ 
of $C^\infty$-functions
 on $\R$
such that%
\footnote{This condition means that functions $\phi$ are smooth as sections
	of line bundles on the projective line $\R\cup\infty$.}
$$
\phi(-1/t) (-t)^{-1+\mu_1-\mu_2\kos  \epsilon_1-\epsilon_2}
$$
also is $C^\infty$-smooth.

Thus, for any fixed $X\in\GL_2(\R)$ we get an operator-valued function
$(\mu_1,\epsilon_1;\mu_2,\epsilon_2)
\mapsto T_{\mu_1,\epsilon_1;\mu_2,\epsilon_2}$
 on $\Lambda$ holomorphic in the variables
$\mu_1$, $\mu_2$.

Generators of the Lie algebra $\frg\frl_2(\R)$ act by formulas
\begin{align}
L_{11}&=-t\frac d{dt}+(-1/2+\mu_1),\qquad
&L_{12}&=\frac d{dt},
\label{eq:lie-R-1}
\\
L_{21}&=-t^2\frac d{dt}+t(-1+\mu_1-\mu_2),
\qquad
&L_{22}&=t\frac d{dt}+(1/2+\mu_2).
\label{eq:lie-R-2}
\end{align}
The expressions for the generators do not depend
on $\epsilon_1$, $\epsilon_2$. However the space of the representation
depends on $\epsilon_1-\epsilon_2$.

Consider integral operators
$$A_{\mu_1,\epsilon_1;\mu_2,\epsilon_2}:
C^\infty_{\mu_1-\mu_2,\epsilon_1-\epsilon_2}
\to
C^\infty_{\mu_2-\mu_1,\epsilon_1-\epsilon_2}
$$
defined by
$$
A_{\mu_1,\epsilon_1;\mu_2,\epsilon_2} f(t):=
\int_{\R} (t-s)^{-1-\mu_1+\mu_2\kos  \epsilon_1-\epsilon_2}
f(s)\,ds.
$$

The integral is convergent if $\Re(\mu_1-\mu_2)<0$ and determines a function
holomorphic in $\mu_1$, $\mu_2$. As usual (see, e.g., \cite{GSh}, \S I.3), the integral admits a meromorphic
continuation to the whole plane $(\mu_1,\mu_2)\in \C^2$ with poles at
$\mu_1-\mu_2=0$, $1$, $2$, \dots. The operators $A_{\mu_1,\epsilon_1;\mu_2,\epsilon_2}$
are intertwinig,
\begin{equation}
A_{\mu_1,\epsilon_1;\mu_2,\epsilon_2}T_{\mu_1,\epsilon_1;\mu_2,\epsilon_2}
=
T_{\mu_2,\epsilon_2;\mu_1,\epsilon_1}A_{\mu_1,\epsilon_1;\mu_2,\epsilon_2}
\label{eq:symmetry}
\end{equation}
(parameters $(\mu_1,\epsilon_1)$ and $(\mu_2,\epsilon_2)$ of $T$ are transposed).
If $\mu_1-\mu_2\ne \pm 1$, $\pm 2$, \dots, then representations $T_{\mu_1,\epsilon_1;\mu_2,\epsilon_2}$
are irreducible and operators $T_{\mu_1,\epsilon_1;\mu_2,\epsilon_2}$ are 
invertible. 

Representations $T_{\mu_1,\epsilon_1;\mu_2,\epsilon_2}$ form a so-called
 (nonunitary) {\it principal series of representations.}  Recall description
of some unitary irreducible representations of $\GL_2(\R)$.

\sm

{\it Unitary principal series.}  If $\mu_1$, $\mu_2\in i\R$, then the representation is unitary 
in $L^2(\R)$.  
Denote by $\Lambda_{\principal}\subset\Lambda$ the subset 
consisting of tuples $(is_1,\epsilon_1;is_2,\epsilon_2)$ with $s_1\ge s_2$.

\sm

{\it Discrete series.} 
Notice that the group $\GL_2(\R)$ acts on the Riemann sphere $\ov C=\C\cup \infty$
by linear-fractional transformations 
\begin{equation}
X:\,z\mapsto \frac{x_{12}+z x_{22}}{x_{11}+z x_{21}},
\label{eq:lin-frac}
\end{equation}
these transformations preserve the real projective line $\R\cup\infty$
and therefore preserves its complement $\C\setminus\R$.
If $\det X>0$, then  (\ref{eq:lin-frac}) leave upper and lower half-planes invariant,
if $\det X<0$, then (\ref{eq:lin-frac}) permutes the half-planes.

Let $n=1$, 2, 3, \dots. Consider the space $H_n$ of 
holomorphic functions $\phi$ on $\C\setminus\R$ satisfying 
$$
\int_{\C\setminus\R} |\phi(z)|^2 |\Im z|^{-n-1} \,d\Re z\,d\Im z<\infty.
$$
In fact, $\phi$ is a pair of holomorphic functions $\phi_+$ and $\phi_-$ determined
  on half-planes $\Im z>0$ and $\Im z<0$. These functions have 
  boundary values on $\R$ in distributional sense (we omit a precise discussion
  since it is not necessary for our purposes).
The space $H_n$ is a Hilbert space with respect to the inner product
$$
\la \phi_1, \phi_2\ra=\int_{\C\setminus\R} \phi_1(z)\ov{\phi_2(z)} \,|\Im z|^{-n-1} \,d\Re z\,d\Im z.
$$
For $s\in\R$, $\delta\in\Z_2$ we define a unitary representation $D_{n,s}$
of $\GL_2(\R)$ in $H_n$ by
\begin{multline*}
D_{n,s}\begin{pmatrix}
x_{11}&x_{12}\\
x_{21}&x_{22}
\end{pmatrix}\phi(z)= \phi\Bigl(\frac{x_{12}+z x_{22}}{x_{11}+z x_{21}}\Bigr)
(x_{11}+z x_{21})^{-1-n}
\times\\\times
\det \begin{pmatrix}
x_{11}&x_{12}\\
x_{21}&x_{22}
\end{pmatrix}^{1/2+n/2+is\kos  \delta}
.
\end{multline*}
In fact, we have operators (\ref{eq:T-mu}) for
$$
\mu_1=-n/2-is,\, \epsilon_1=-1+n+\delta;\,\, \mu_2=n/2+is,\, \epsilon_2=\delta
$$
restricted to the subspace generated by  boundary values of 
functions $f_+$ and $f_-$.

We denote by $\Lambda_{\discrete}$ the set of all parameters $(n,is)$ 
of discrete series. 

\sm




\sm

{\bf \punct The Fourier transform.%
	\label{ss:R-Fourier}} Let $\phi$ be contained
in the space $C_0^\infty(\GL_2(\R)$ of compactly supported function on
$\GL_2(\R)$.  We consider
a function that 
sent each $\wt \mu=(\mu_1,\epsilon_1;\mu_2,\epsilon_2)$ to
an operator in $C^\infty_{\mu_1-\mu_2,\epsilon_1-\epsilon_2}$
given by
$$
T_{\mu_1,\epsilon_1;\mu_2,\epsilon_2}(F)=
\int_{\GL_2(\R)} F(X)\, T_{\mu_1,\epsilon_1;\mu_2,\epsilon_2}(X)\,\frac {dX}{\det(X)^2}
.
$$

Next, we define a subset $\Lambda_{\tempered}\subset \Lambda$ by $$\Lambda_{\tempered}:=\Lambda_{\principal}\cup\Lambda_{\discrete}.$$
Let us define {\it Plancherel measure} $d\cP(\wt\mu)$ on 
$\Lambda_{\tempered}$.
On the piece $\Lambda_{\principal}$ it is given by
\begin{align*}
d\cP(is_1,0;is_2,\epsilon_2)
=\frac 1{16\pi^3}(s_1-s_2)
\tanh \pi(s_1-s_2)/2\,ds_1\,ds_2;
\\
d\cP(is_1,1;is_2,\epsilon_2)
= \frac 1{16\pi^3}
(s_1-s_2)
\coth \pi(s_1-s_2)/2\,ds_1\,ds_2.
\end{align*}
On $n$-th piece of $\Lambda_{\discrete}$ it is given by
$$
d\cP=\frac n{8\pi^3} ds.
$$

Consider the space $\cL^2$ of
 functions $Q$ on $\Lambda_{\tempered}$ taking values
in the space of Hilbert--Schmidt operators in the corresponding Hilbert spaces
and
satisfying the condition
$$
\int_{\Lambda_{\tempered}} \tr\Bigl( Q(\wt\mu)Q^*(\wt\mu)\Bigr)\,d\cP(\wt\mu)<\infty.
$$
This is a Hilbert space with respect to the inner product
$$
\la Q_1,Q_2\ra_{\cL^2}:= \int_{\Lambda_{\tempered}}
 \tr \Bigl(Q_1(\wt\mu)Q_2^*(\wt\mu)\Bigr)\,\,d\cP(\wt\mu).
$$
According the Plancherel theorem, {\it for any $F_1$, $F_2\in C^\infty_0(\GL_2(\R))$
we have 
$$
\la F_1, F_2\ra_{L^2(\GL_2(\R))}
=\la  T(F_1), T(F_2)\ra_{\cL^2(\Lambda_{\tempered})}.
$$
Moreover, the map $F\to T(F)$ extends to a unitary operator
$L^2(\GL_2(\R))\to \cL$.}

\sm

{\bf\punct Overgroup.%
	\label{ss:R-overgroup}}  Let $\Mat_2(\R)$ be the space of all real
matrices of order 2.
By $\Gr_{4,2}(\R)$ we denote the Grassmannian of 2-dimensional subspaces
in $\R^2\oplus \R^2$. For any operator $\R^2\to\R^2$ 
 its graph is an element $\Gr_{4,2}(\R)$. The set $\Mat_2(\R)$ of such operators
is an open dense chart in $\Gr_{4,2}(\R)$.

 The group $\GL_4(\R)$ acts in a natural way in $\R^4$ and therefore 
on the Grassmannian. In the chart $\Mat_2(\R)$
the action is given by the formula (see, e.g., \cite{Ner-gauss}, Theorem 2.3.2)
$$
\begin{pmatrix} A&B\\C&D \end{pmatrix}:\,X\mapsto (A+XC)^{-1}(B+XD),
$$
where $\begin{pmatrix} A&B\\C&D \end{pmatrix}$
is an element of $\GL_4(\R)$ written as a  block matrix of size
 $2+2$. The Jacobian of this transformation
 is (see, e.g., \cite{Ner-gauss},  Theorem 2.3.2)
 $$
 \det(A+XC)^{-4}\det \begin{pmatrix} A&B\\C&D \end{pmatrix}^2.
 $$
 For $\sigma\in i\R$ we define a unitary representation
 of $\GL_4(\R)$ in $L^2(\Mat_2(\R))$ by
 \begin{multline*}
 R_\sigma \begin{pmatrix} A&B\\C&D \end{pmatrix} F(X)=F\bigl( (A+XC)^{-1}(B+XD)\bigr)
 \times\\\times
 \,|\det(A+XC)|^{-2+2\sigma} 
 \left| 	\det\begin{pmatrix} A&B\\C&D \end{pmatrix} \right|^{1-\sigma}
 \end{multline*}
 (these representations are contained in degenerate principal series).
 
 The group $\GL_2(\R)$ is an open dense subset in $\Mat_2(\R)$.
 Therefore, we can identify the spaces $L^2$ on $\GL_2(\R)$ 
 and $\Mat_2(\R)$. For this we consider a unitary operator
 $$
 J:\, L^2\bigl(\Mat_2(\R)\bigr)\to L^2\bigl(\GL_2(\R) \bigr)
 $$
 given by
 $$
 J_\sigma F(X)=F(X)\cdot |\det X|^{1-\sigma}
 $$
This determines a unitary representation 
$U_\sigma:=J_\sigma R_\sigma J_\sigma^{-1}$
of $\GL_4(\R)$ in $L^2\bigl(\GL_2(\R) \bigr) $,
the explicit expression is 
\begin{multline}
U_\sigma \begin{pmatrix} A&B\\C&D \end{pmatrix} 
F(X)=F\bigl( (A+XC)^{-1}(B+XD)\bigr) 
\times\\\times
\Biggl|\frac{\det(A+XC)\det(B+XD)}{\det X 
	\det\begin{pmatrix} A&B\\C&D \end{pmatrix}}  \Biggr|^{-1+\sigma}
.
\label{eq:U-sigma}
\end{multline}

Consider a subgroup $\GL_2(\R)\times \GL_2(\R)\subset \GL_4(\R)$
consisting of matrices $\begin{pmatrix} A&0\\0&D \end{pmatrix} $.
For this subgroup we get the usual left-right action in 
$L^2\bigl(\GL_2(\R)\bigr)$,
 $$
 U_\sigma \begin{pmatrix} A&0\\0&D \end{pmatrix}
 F(X)=F(A^{-1}X D).
 $$
 Formula (\ref{eq:U-sigma}) extends this formula to the whole group
 $\GL_4(\R)$. The Lie algebra $\frg\frl_4$ acts in the space
 of functions on $\Mat_2(\R)$ by first order differential operators,
 which can be easily written; a list of formulas for all
 generators $e_{kl}$, where $k$, $l=1$, 2, 3, 4, is given below in Subs. (\ref{ss:list}).
 We restrict this action to the space of smooth compactly supported functions
 on $\GL_2(\R)$. Notice that the operators $i\cdot e_{kl}$ are symmetric on this domain,
 but some of them are not essentially self-adjoint.

Our purpose is to write explicitly the images $E_{kl}$ of operators $e_{kl}$
 under the Fourier
transform. 
 
 \sm
 
 {\bf \punct Formulas.%
 	\label{ss:R-formulas}}
 Operators $T_{\mu_1,\epsilon_1;\mu_2;\epsilon_2}(F)$
 have the form
 \begin{equation}
 T_{\mu_1,\epsilon_1;\mu_2;\epsilon_2}(F)\phi(t)=
 \int_{-\infty}^{\infty} K(t,s|\mu_1,\epsilon_1;\mu_2,\epsilon_2)\,
 \phi(t)\,dt.
 \end{equation}
 Recall that functions $K$ are holomorphic in $\mu_1$, $\mu_2$.

 We wish to write operators $E_{kl}$ on kernels $K$.
 The complete list is contained below in Subsect. \ref{ss:list},
 here we present two basic expressions.
 
 The algebra $\frgl_{4}(\R)$ can be decomposed as a linear space
 into a direct sum of four subalgebras $\fra$, $\frb$, $\frc$, $\frd$
 consisting of matrices of the form 
 $$
 \begin{pmatrix}
 *&0\\0&0
 \end{pmatrix},
 \qquad
  \begin{pmatrix}
 0&0\\0&*
 \end{pmatrix},
 \qquad \begin{pmatrix}
 0&0\\ *&0
 \end{pmatrix},
 \qquad \begin{pmatrix}
 0&0\\0&*
 \end{pmatrix}
 $$ 
 The subalgebras $\fra$ and $\frd$ are isomorphic to $\frgl_2(\R)$,
  subalgebras $\frb$ and $\frc$ are Abelian.

  Formulas for the action of $\fra$ and $\frd$
    immediately follow from the definition of the Fourier transform.
    To obtain formulas for the whole $\frgl_4$, it is sufficient
    to write expressions for one generator of $\frb$, say $E_{14}$,
     and one generator
    of $\frc$, say $E_{32}$. After this
     other generators can be obtained by evaluation of commutators.

 We define shift operators $V_1^+$, $V_1^-$, $V_2^+$, $V_2^-$
 by
 \begin{align*}
 V_1^\pm  K(t,s|\mu_1,\epsilon_1;\mu_2,\epsilon_2)= K(t,s|\mu_1\pm 1,\epsilon_1+1;\mu_2,\epsilon_2);
 \\
 V_2^\pm  K(t,s|\mu_1,\epsilon_1;\mu_2,\epsilon_2)= K(t,s|\mu_1,\epsilon_1;\mu_2\pm 1,\epsilon_2+1).
 \end{align*}
 In this notation,
  \begin{align*}
  E_{14}=\frac{-1/2 - \sigma + \mu_1}{\mu_1 - \mu_2}\, \frac\partial{\partial s}\, V_1^- 
  +\frac{-1/2 - \sigma + \mu_2}{\mu_1 - \mu_2}\,\frac\partial{\partial t}\, V_2^-;
  \end{align*}
 
     \begin{equation*}
     E_{32}=\frac{1/2 +  \mu_1 + \sigma}{\mu_1 - \mu_2}\,
     \frac\partial{\partial t}\, V_1^+  +\frac{1/2 +  \mu_2 +  \sigma}{\mu_1 - \mu_2} \, \frac\partial{\partial s}\,
     V_2^+.
     \end{equation*} 


     {\bf \punct Remarks on a general problem.} 
    In \cite{Ner-over} the author formulated the following question:
           {\it  Assume that we know the explicit Plancherel formula for the
           	restriction of a unitary representation $\rho$ of a group G to a subgroup H. Is
           	it possible to write the action of the Lie algebra of G in the direct integral of
           	representations of H?} 
           
           Now it seems that an answer to this question is affirmative.

           The initial paper \cite{Ner-over} contains a solution for a tensor product%
           \footnote{The problem of a decomposition  of  a tensor product
           	$\rho_1\otimes\rho_2$
           of representations of a group $G$ is a special case of a decomposition of restrictions.
            Indeed,  $\rho_1\otimes\rho_2$ is a representation of the group $G\times G$, and
            we must restrict this representation to the diagonal $G$.} of a highest and lowest weight representations of $\SL_2(\R)$. In this case the  overalgebra
           acts by differential-difference operators in the space $L^2(S^1\times \R_+)$
            having the form
           $$
           Lf(\phi,s)=D_1 f(\phi, s+i)+ D_2 f(\phi,s) +D_3 f(\phi,s-i),
           $$
           where $D_j$ are second order differential operators in the variable $\phi$.

           In \cite{Mol1}--\cite{Mol4} Molchanov solved several rank 1 problems of this type,
 expressions are similar, but there appear differential operators of order 4.
 In \cite{Ner-gl} the action of the overalgebra in restrictions
 from $\GL_{n+1}(\C)$ to $\GL_n(\C)$, in this case differential operators have order $n$.  
 In all the cases examined by now shift operators in the imaginary direction are present.   
           
 In the present paper, we write the action of the overalgebra    
 in the restriction of a degenerate principal series of the group
 $\GL_4(\R)$ to $\GL_2(\R)$. Notice that canonical overgroups exist for all
 10 series of real classical groups%
 \footnote{Emphasize that in the present paper
 	 we consider groups $\GL$ and not $\SL$, also 
 	we must consider groups $\U(p,q)$ and not $\SU(p,q)$.}.
 Moreover overgroups exist for all 52 series of classical semisimple symmetric
 spaces $G/H$, see \cite{Mak}, \cite{Ner-uniform}, see also
 \cite{Ner-gauss}, Addendum D.6. So the problem makes sense for all
 classical symmetric spaces.

 Shturm--Liouville problems for difference operators in $L^2(\R)$ in
 the imaginary direction
 $$
 \lambda f(s)=a(s) f(s+i)+b(s) f(s)+ c(s) f(s-i)
 $$
  arise in a natural way  in the theory of hypergeometric orthogonal polynomials, see, e.g.,
  \cite{AAR}, \cite{KSh}, apparently a first example (the Meixner--Pollaszek system)
   was discovered by J.~Meixner
  in 1930s. On such operators with continuous spectra see \cite{Ner-index},
  \cite{Gro},
  \cite{Ner-imaginary}. See also a multi-rank work of I.~Cherednik 
  \cite{Cher}
  on Harish-Chandra spherical transforms.
  
  \sm
  
    {\bf\punct The Fourier transform on $\GL_2(\C)$.%
    	\label{ss:C-Fourier}} For a detailed exposition 
    of representations of the Lorentz group, see \cite{Nai}.
  For $\nu$, $\nu'\in\C$ satisfying $\nu-\nu'\in\Z$ we define the function
  $z^{\nu\|\nu'}$ on the multiplicative group of $\C$ by
  $$
  z^{\nu\|\nu'}:=z^\nu z^{\nu'}:=|z|^{2\nu}\, \ov z^{\,\nu'-\nu}.
  $$
Denote by $\Delta$ the set of all
$
(\mu_1,\mu_1';\mu_2,\mu_2')\in\C^4$
such that $\mu_1-\mu_1'\in\Z$, $\mu_2-\mu_2'\in\Z$.
For $(\mu_1,\mu_1';\mu_2,\mu_2')\in \Delta$
we define a representation of $\GL_2(\C)$ in the space of functions on
$\C$ by 
\begin{multline*}
T_{\mu_1,\mu_1';\mu_2,\mu_2'}
\begin{pmatrix}
x_{11}&x_{12}\\
x_{21}&x_{22}
\end{pmatrix} \phi(t)=\\=
\phi\Bigl(\frac{x_{12}+t x_{22}}{x_{11}+t x_{21}}\Bigr)
(x_{11}+t x_{21})^{-1+\mu_1-\mu_2\|-1+\mu_1'-\mu_2'} \det\begin{pmatrix}
x_{11}&x_{12}\\
x_{21}&x_{22}
\end{pmatrix}^{1/2+\mu_2\|1/2+\mu_2'}
.
\end{multline*}
We consider the space of $C^\infty$-smooth functions $\phi$
on $\C$ such that
$$
\phi(-t^{-1}) t^{-1+\mu_1-\mu_2\|-1+\mu'_1-\mu'_2}
$$
is $C^\infty$-smooth at 0.
These representations form the (nonunitary) {\it principal series}.

It is convenient to complexify the Lie algebra of $\GL_2(\C)$,
$$
\frgl_2(\C)_\C\simeq \frgl_2(\C)\oplus \frgl_2(\C).
$$
Under this isomorphism, the operators of the Lie algebra act in our representation
by
\begin{align}
L_{11}&=-t\frac d{dt}+(-1/2+\mu_1),\qquad
&L_{12}&=\frac d{dt},
\label{eq:lie-C-1}
\\
L_{21}&=-t^2\frac d{dt}+t(-1+\mu_1-\mu_2),
\qquad
&L_{22}&=t\frac d{dt}+(1/2+\mu_2),
\label{eq:lie-C-2}
\\
\ov L_{11}&=-\ov t\frac d{ d\ov t}+(-1/2+\mu_1'),\qquad
&\ov L_{12}&=\frac d{d\ov t},
\label{eq:lie-C-3}
\\
\ov L_{21}&=-\ov t^2\frac d{d\ov t}+\ov t(-1+\mu_1'-\mu_2'),
\qquad
&\ov L_{22}&=\ov t\frac d{d\ov t}+(1/2+\mu_2').
\label{eq:lie-C-4}
\end{align}
Formally, we have duplicated expressions (\ref{eq:lie-R-1})--(\ref{eq:lie-R-2}).
     
If 
\begin{equation}
\Re (\mu_1+\mu_1')=0,\qquad\Re(\mu_2+\mu_2')=0,
\label{eq:c-tempered}
\end{equation}
 then  the representation
$T_{\mu_1,\mu_1';\mu_2,\mu_2'}$ is unitary in $L^2$.
Denote by $\Delta_{\tempered}$ the set
of such tuples $(\mu_1,\mu_1';\mu_2,\mu_2')$, it is a union of
a countable family of parallel 2-dimensional real planes in $\C^4$,
we equip it by a natural Lebesgue measure $d\lambda(\mu)$.

For any compactly supported smooth function $F$ on
$\GL_2(\C)$ we define its {\it Fourier transform}
as an operator-valued function on $\Delta$ given by 
$$
T_{\mu_1,\mu_1';\mu_2,\mu_2'} (F) =\int_{\GL_2(\C)}
F(X)\, |\det X|^{-4}\,  \prod_{k,l=1,2} d\Re x_{kl}\,d\Im x_{kl}.
$$

The Plancherel formula is the following identity
\begin{multline*}
\la F_1,F_2\ra_{L^2\bigl(\GL_2(\C)\bigr)}
=-C\cdot\int_{\Delta_{\tempered}} \tr \bigl( T_{\mu_1,\mu_1';\mu_2,\mu_2'} (F_1)
T_{\mu_1,\mu_1';\mu_2,\mu_2'} (F_2)\bigr)^*
\times\\\times (\mu_1-\mu_2)  (\mu_1'-\mu_2')
d\lambda(\mu),
\end{multline*}
where $C$ is an explicit constant.

Denote by $K(t,s|\mu_1,\mu_1';\mu_2,\mu_2')$ the kernel of the operator
$T_{\mu_1,\mu_1';\mu_2,\mu_2'} (F)$,
$$
T_{\mu_1,\mu_1';\mu_2,\mu_2'} (F) \phi(t)=
\int_\C  K(t,s|\mu_1,\mu_1';\mu_2,\mu_2')\phi(s)\,ds.
$$

{\bf \punct Overgroup for $\GL_2(\C)$.%
	\label{ss:C-overgroup}}
Consider the complex Grassmannian $\Gr_{4,2}(\C)$ of 2-dimensional
planes in $\C^4$, again the set $\Mat_2(\C)$ is an open dense set
on $\Gr_{4,2}(\C)$.
For $\sigma$, $\sigma'\in i \R$ consider a unitary representation
$R_{\sigma,\sigma'}$
of $\GL_4(\C)$ in $L^2\bigl(\Mat_2(\C)\bigr)$ 
given by 
\begin{multline*}
R_{\sigma,\sigma'} F(X)
=\begin{pmatrix} A&B\\C&D \end{pmatrix} F(X)=F\bigl( (A+XC)^{-1}(B+XD)\bigr)
\times\\\times
\,|\det(A+XC)|^{-2+2\sigma\|-2+2\sigma'} 
\left| 	\det\begin{pmatrix} A&B\\C&D \end{pmatrix} \right|^{1-\sigma\|1-\sigma'}
.
\end{multline*}
We define a unitary operator $J:L^2\bigl(\Mat_2(\C)\bigr)
\to L^2\bigl(\Mat_2(\C)\bigr)$ by
$$
J_{\sigma,\sigma'} F(X)=F(X)(\det X)^{-1+\sigma\|-1+\sigma'}.
$$
In this way we get a unitary representation
$
U_{\sigma,\sigma'}=J_{\sigma,\sigma'}  R_{\sigma,\sigma'} J_{\sigma,\sigma'}^{-1} 
$
of $\GL_4(\C)$ in $L^2\bigl(\Mat_2(\C)\bigr)$:
 \begin{multline*}
U_{\sigma\|\sigma'} \begin{pmatrix} A&B\\C&D \end{pmatrix} 
F(X)=F\bigl( (A+XC)^{-1}(B+XD)\bigr) 
\times\\\times
\Biggl|\frac{\det(A+XC)\det(B+XD)}{\det X 
	\det\begin{pmatrix} A&B\\C&D \end{pmatrix}}  \Biggr|^{-1+\sigma\|-1+\sigma'}
.
 \end{multline*}
 
 {\bf \punct Formulas for $\GL_2(\C)$.%
 	\label{ss:C-formulas}}
We wish to write the action of the Lie algebra 
$$
\frgl_4(\C)_\C\simeq \frgl_4(\C)\oplus \frgl_4(\C)
$$
in the Plancherel decomposition of $\GL_2(\C)$.
Denote the standard generators of $\frgl_4(\C)\oplus 0$
and $0\oplus\frgl_4(\C)$  by $E_{kl}$ and $\ov E_{kl}$
respectively.
Define the following shift operators
\begin{align*}
V_1 K(t,s;\mu_1,\mu_1',\mu_2,\mu_2')= K(t,s;\mu_1+1,\mu_1',\mu_2,\mu_2');
\\
V_1' K(t,s;\mu_1,\mu_1',\mu_2,\mu_2')= K(t,s;\mu_1,\mu_1'+1,\mu_2,\mu_2'),
\end{align*}
and similar operators $V_2$ and $V_2'$ shifting $\mu_2$ and $\mu_2'$.

 Then
   \begin{align*}
   E_{14}&=\frac{-1/2 - \sigma + \mu_1}{\mu_1 - \mu_2}\, \frac\partial{\partial s} \, V_1^{-1} 
   -\frac{-1/2 - \sigma + \mu_2}{\mu_1 - \mu_2}\,\frac\partial{\partial t}\, V_2^{-1}
   ;
   \\
     \ov E_{14}&=\frac{-1/2 - \sigma' + \mu'_1}{\mu'_1 - \mu'_2}
     \, \frac\partial{\partial \ov s} \,(V_1')^{-1} 
     -\frac{-1/2 - \sigma' + \mu'_2}{\mu'_1 - \mu'_2}\,
     \frac\partial{\partial \ov t}\, (V_2')^{-1}
    ; \\
   E_{32}&=\frac{1/2 +  \mu_1 + \sigma}{\mu_1 - \mu_2}
  \, \frac\partial{\partial t}\, V_1  +\frac{1/2 +  \mu_2 +  \sigma}{\mu_1 - \mu_2}  \,\frac\partial{\partial s}
  \, V_2
;\\
  \ov  E_{32}&=\frac{1/2 +  \mu_1' + \sigma'}{\mu'_1 - \mu'_2}\,
    \frac\partial{\partial \ov t}\, V_1'  +\frac{1/2 +  \mu'_2 +  \sigma'}{\mu'_1 - \mu'_2} \, \frac\partial{\partial \ov s}\,
    V_2'
 .   \end{align*}  
 
 \section{Calculations}
 
 \COUNTERS
 
 {\bf\punct The expression for kernel.}
 
 \begin{lemma}
 	The kernel $K(\cdot)$ of an integral operator
 	$T_{\mu_1,\epsilon_1;\mu_2,\epsilon_2}(F)$
 	is given by the
 	formula
 	\begin{multline}
 	 K(t,s|\mu_1,\epsilon_1;\mu_2,\epsilon_2)=\\=
 	 \iiint\limits_{\R^3} F\bigl(
 	u - t v, s u - s t v - t w, v, s v + w \bigr)\,u^{-3/2+\mu_1\kos  \epsilon_1}
 	w^{-3/2+\mu_2\kos \epsilon_2}\,du\,dv\,dw.
 	\label{eq:K}
 	 \end{multline}
 	 For $F\in C_0^\infty\bigl(\GL_2(\R)\bigr)$
 	 the integration is actually taken over a bounded domain.
 \end{lemma}
 
 The integral converges if $\Re \mu_1>1/2$, $\Re \mu_2>1/2$.
 For fixed $\epsilon_1$, $\epsilon_2$
 it has a meromorphic continuation to the whole
 complex plane $(\mu_1,\mu_2)$ with poles on the hyperplanes $\mu_1=-1/2-k$,
 $\mu_2=-1/2-k$, where $k=0$, 1, 2, \dots (see, e.g., \cite{GSh}, \S I.3).
 
 \sm

 {\sc Proof.} By the definition
 \begin{multline*}
 \int\limits_{\R} 
 T_{\mu_1,\epsilon_1;\mu_2,\epsilon_2}(F)\phi(t)\,\,\psi(t)\,dt
 =\\=
 \int\limits_{\R}\int\limits_{\Mat_2(\R)}
 F(x_{11},x_{12},x_{21},x_{22})\phi\Bigl(\frac{x_{12}+t x_{22}}{x_{21}+t x_{12}}\Bigr)
 (x_{11}+t x_{21})^{-1+\mu_1-\mu_2\kos \epsilon_1-\epsilon_2}
 \times\\\times
  (x_{11}x_{22}-x_{12}x_{21})^{1/2+\mu_2\kos \epsilon_2}
  \,\frac{dx_{11}\,dx_{12}\,dx_{21}\,dx_{22}}{(x_{11}x_{22}-x_{12}x_{21})^2}\,dt.
 \end{multline*}
 
 In the interior integral, we pass from the variables $x_{11}$, $x_{12}$, $x_{21}$,
 $x_{22}$ to new variables $u$, $v$, $w$, $s$ defined by
 \begin{equation}
 \begin{pmatrix}
 x_{11}&x_{12}\\x_{21}&x_{22}
 \end{pmatrix}
 =
 \begin{pmatrix}
 1&-t\\0&1
 \end{pmatrix}
 \begin{pmatrix}
 u&0\\v&w
 \end{pmatrix}
  \begin{pmatrix}
  1&s\\0&1
  \end{pmatrix}
  \label{eq:zamena}
 \end{equation}
 or
 $$
 x_{11}=u - t v,\qquad x_{12}= s u - s t v - t w,\qquad x_{21}= v,\qquad x_{22}= s v + w
 .
 $$
  The Jacobi matrix of this transformation is triangular, and the Jacobian
 is $|u|$. The inverse transformation is
 $$
 u=x_{11}+t x_{21},\quad v=x_{21}, \quad w=\frac{x_{11}x_{22}-x_{12}x_{21}}{x_{11}+t x_{21}},
 \qquad
s= \frac{x_{12}+t x_{22}}{x_{11}+t x_{21}}
. $$
 We also have $x_{11}x_{22}-x_{12}x_{21}=uw$.
 
 After the change of variables we come to
 $$
 \iint_{\R^2} K(t,s|\mu_1,\epsilon_1;\mu_2,\epsilon_2)\,
 \phi(s)\,\psi(t)\,ds\,dt,
 $$
 where $K(\cdot)$ is given by (\ref{eq:K}).
 
 A function $F$ has a compact support in 
 $\R^4\setminus\{x_{11}x_{22}-x_{12}x_{21}=0\}$.
 So, actually, $x_{21}=v$, $x_{11}=u-tv$, $x_{22}=w+sv$
 are contained in a bounded domain. This implies the second claim
 of the lemma.
 \hfill $\square$
 
 \sm
 
 {\bf\punct Preliminary remarks.}
 Below $F$ denotes
 $$
 F:=F(u - t v, s u - s t v - t w, v, s v + w).
 $$
Also, $\partial_{11}F$, $\partial_{12}F$, etc. denote
$$
\partial_{11}F:=\frac{\partial}{\partial x_{11}}
F(x_{11}, x_{12}, x_{21}, x_{22})\Biggr|_{\begin{matrix}
x_{11}= u - t v& x_{12}=s u - s t v - t w ,\\
 x_{21}=v,&
 x_{22}=s v + w	
	\end{matrix}
	}
$$
etc.
Partial derivatives of $F$ are
\begin{align}
\frac{\partial}{\partial u} F&=\partial_{11} F+s\, \partial_{12} F;
\label{eq:partial-u}
\\
\frac{\partial}{\partial v} F&=-t\, \partial_{11} F-st\, \partial_{12} F+
\partial_{21} F+s\, \partial_{22} F;
\label{eq:partial-v}
\\
\frac{\partial}{\partial w} F&=-t\,\partial_{12}F+ \partial_{22} F,
\label{eq:partial-w}
\end{align}
 and
 \begin{align}
 \frac{\partial}{\partial s} F& = (u-vt)\, \partial_{12}F+v\, \partial_{22}F;
 \label{eq:partial-s}
 \\
 \frac{\partial}{\partial t} F&=-v\, \partial_{11}F-(w+sv)\,\partial_{12}F.
 \label{eq:partial-t}
 \end{align}
 Also, notice that
 $$
 \left(y^{\nu\kos \delta}\right)'=\nu y^{\nu-1\kos \delta+1} .
 $$

 \sm
 
 {\bf\punct A verification of the formula for
 	$E_{14}$.} It is easy to verify that
 the operator $e_{14}$ in $C^\infty_0\bigl(\GL_2(\R)\bigr)$
 is given by
 $$
 e_{14}= \frac \partial{\partial x_{12}} - (-1 + \sigma)\frac{x_{21}}{x_{11}x_{22}-x_{12}x_{21}}
. $$
 Therefore,
 \begin{multline}
 T_{\mu_1,\epsilon_1;\mu_2,\epsilon_2}(e_{14}F)=\\=
  	 \iiint \partial_{12} F
  	\, \,u^{-3/2+\mu_1\kos \epsilon_1}
  	 w^{-3/2+\mu_2\kos \epsilon_2}\,du\,dv\,dw-\\-
  (-1 + \sigma)	 
    	 \iiint  F
    	  \frac v{uw} \,u^{-3/2+\mu_1\kos \epsilon_1}
    	 w^{-3/2+\mu_2\kos \epsilon_2}\,du\,dv\,dw,
    	 \label{eq:T14}
 \end{multline}
(the integration is taken over $\R^3$ on default).
 
 We must verify that (\ref{eq:T14}) coincides with
 \begin{equation}
 E_{14} K= \Bigl(\frac{-1/2 - \sigma + \mu_1}{\mu_1 - \mu_2}\, \frac\partial{\partial s}
 \, V_1^- 
 +\frac{-1/2 - \sigma + \mu_2}{\mu_1 - \mu_2}\,\frac\partial{\partial t}\, V_2^-\Bigr)
 K
 \label{eq:T14-1}
 \end{equation}
 
Below we establish two formulas
   	\begin{multline}
   	 \iiint \partial_{12} F\, \,u^{-3/2+\mu_1\kos \epsilon_1}
   	 w^{-3/2+\mu_2\kos \epsilon_2}\,du\,dv\,dw
\,\, -\,\,	\frac\partial{\partial s} V_1^- K
 =\\=
(-3/2+\mu_2)  \iiint  F \,\,
  \frac v{uw} \,u^{-3/2+\mu_1\kos \epsilon_1}
 w^{-3/2+\mu_2\kos \epsilon_2}\,du\,dv\,dw,
 \label{eq:1}
  \end{multline}
  \begin{multline}
     	 \iiint \partial_{12} F\,\,u^{-3/2+\mu_1\kos \epsilon_1}
     	 w^{-3/2+\mu_2\kos \epsilon_2}\,du\,dv\,dw
     	 \,\, +\,\,	\frac\partial{\partial t} V_2^- K
     	  =\\=
     	  (-3/2+\mu_1)  \iiint F\,\,  \frac v{uw} \,u^{-3/2+\mu_1\kos \epsilon_1}
     	  w^{-3/2+\mu_2\kos \epsilon_2}\,du\,dv\,dw.
     	  \label{eq:2}
  \end{multline}
 Considering the sum of (\ref{eq:1}) and (\ref{eq:2})
 with coefficients
 $\frac{-1/2 - \sigma + \mu_1}{\mu_1 - \mu_2}$ and
 $\frac{-1/2 - \sigma + \mu_2}{\mu_1 - \mu_2}$
 we get coincidence of (\ref{eq:T14}) and (\ref{eq:T14-1});
 for this, we use he identities
 \begin{align*}
& \frac{-1/2 - \sigma + \mu_1}{\mu_1 - \mu_2}-\frac{-1/2 - \sigma + \mu_2}{\mu_1 - \mu_2}=1;
 \\
 &(-3/2+\mu_2)\cdot  \frac{-1/2 - \sigma + \mu_1}{\mu_1 - \mu_2}-
 (-3/2+\mu_1)\cdot  \frac{-1/2 - \sigma + \mu_2}{\mu_1 - \mu_2}=1-\sigma.
 \end{align*}
 
 Now let us check (\ref{eq:1}).
 The following identity can be verified by a straightforward calculation
 (with (\ref{eq:partial-s}) and (\ref{eq:partial-w})):
 $$
 \partial_{12}F-\frac 1 u \frac \partial{\partial s} F
 =-\frac v u \frac \partial{\partial w} F.
 $$
 Therefore the left-hand side of (\ref{eq:1}) equals to
 $$
 \iiint \Bigl[ -\frac v u \frac \partial{\partial w} F \Bigr]\cdot 
 \,u^{-3/2+\mu_1\kos \epsilon_1}
 w^{-3/2+\mu_2\kos \epsilon_2}\,du\,dv\,dw.
 $$
 We integrate this expression by parts in the variable $w$ and come to (\ref{eq:1}).
 
 To check (\ref{eq:2}), we verify the identity 
 $$
 \partial_{12} F +\frac 1 v \frac{\partial}{\partial t} F
 =\frac vw \frac{\partial}{\partial u} F
 $$
 and after this integrate by parts as above.
 
 \sm
 
 {\bf \punct  A verification of the formula for $E_{32}$.}
 We have
$$ e_{32} = -\Bigl(x_{11} x_{21} \frac\partial  {\partial x_{11}} + 
 x_{11} x_{22} \frac\partial {\partial x_{12}} + 
 x_{21}^2 \frac\partial {\partial x_{21}} + 
 x_{21} x_{22} \frac\partial {\partial x_{22}}\Bigr) + (-1 + \sigma)x_{21}
 .
 $$
 Therefore,
 \begin{equation}
 T(e_{32}K)=
 \iiint \bigl(G+ (-1+\sigma)v F\bigr) \cdot 
 \,u^{-3/2+\mu_1\kos \epsilon_1}
 w^{-3/2+\mu_2\kos \epsilon_2}\,du\,dv\,dw,
 \label{eq:e32K}
 \end{equation}
 where
 $$
 G=(u-vt) v\,\partial_{11} F+ (u-vt)(w+vs) \,\partial_{12} F +v^2\, \partial_{21} F
 +v (w+vs) \,\partial_{22} F.
 $$
 On the other hand,
     \begin{multline}
     E_{32}K=\frac{1/2 +  \mu_1 + \sigma}{\mu_1 - \mu_2}
     \iiint
     u\frac\partial{\partial t}F
     \cdot 
     \,u^{-3/2+\mu_1\kos \epsilon_1}
     w^{-3/2+\mu_2\kos \epsilon_2}\,du\,dv\,dw 
     +\\  +
     \frac{1/2 +  \mu_2 +  \sigma}{\mu_1 - \mu_2}  
     \iiint w\frac\partial{\partial s} F
     \cdot 
     \,u^{-3/2+\mu_1\kos \epsilon_1}
     w^{-3/2+\mu_2\kos \epsilon_2}\,du\,dv\,dw.
     \label{eq:E32-K}
     \end{multline}
     We must verify that (\ref{eq:e32K}) and (\ref{eq:e32K}) are equal.
 As in the previous subsection, this statement is reduced to a pair of identities
 \begin{multline}
 \iiint \bigl(G+ (-1+\sigma)v F-u\frac{\partial}{\partial t} F\bigr) \cdot 
 \,u^{-3/2+\mu_1\kos \epsilon_1}
 w^{-3/2+\mu_2\kos \epsilon_2}\,du\,dv\,dw
 =\\=
 (1/2+\mu_2+\sigma)
  \iiint v F \cdot 
 \,u^{-3/2+\mu_1\kos \epsilon_1}
 w^{-3/2+\mu_2\kos \epsilon_2}\,du\,dv\,dw;
 \label{eq:3}
 \end{multline}
  \begin{multline}
  \iiint \bigl(G+ (-1+\sigma)v F+w\frac{\partial}{\partial s} F\bigr) \cdot 
  \,u^{-3/2+\mu_1\kos \epsilon_1}
  w^{-3/2+\mu_2\kos \epsilon_2}\,du\,dv\,dw
  =\\=
  (1/2+\mu_1+\sigma)
  \iiint v F \cdot 
  \,u^{-3/2+\mu_1\kos \epsilon_1}
  w^{-3/2+\mu_2\kos \epsilon_2}\,du\,dv\,dw.
  \label{eq:4}
  \end{multline}
 Let us verify (\ref{eq:3}). It can be easily checked (with (\ref{eq:partial-t}),
 (\ref{eq:partial-v}), (\ref{eq:partial-w})) that
$$
 G -u\frac{\partial}{\partial t} F=-v^2 \frac{\partial}{\partial v}F-vw \frac{\partial}{\partial w} F.
$$
We substitute this to the left-hand side of (\ref{eq:3}) and
come to
 \begin{multline}
 \iiint \bigl(-v^2 \frac{\partial}{\partial v}F-vw \frac{\partial}{\partial w} F  + (-1+\sigma)v F\bigr) \times\\\times
 \,u^{-3/2+\mu_1\kos \epsilon_1}
 w^{-3/2+\mu_2\kos \epsilon_2}\,du\,dv\,dw.
 \end{multline}
Integrating by parts, we get
\begin{multline*}
\iiint F
\cdot \frac{\partial}{\partial v} \Bigl(v^2\Bigr)
\cdot u^{-3/2+\mu_1\kos \epsilon_1}
w^{-3/2+\mu_2\kos \epsilon_2}\,du\,dv\,dw
+\\+
\iiint vF
\cdot u^{-3/2+\mu_1\kos \epsilon_1}
\frac{\partial}{\partial v} \Bigl(w^{-1/2+\mu_2\kos \epsilon_2+1}\Bigr)\,du\,dv\,dw
+\\
+
(-1+\sigma) \iiint  v F \cdot 
 \,u^{-3/2+\mu_1\kos \epsilon_1}
 w^{-3/2+\mu_2\kos \epsilon_2}\,du\,dv\,dw.
\end{multline*}
After a summation we come to the right-hand side of (\ref{eq:3}).

A proof of (\ref{eq:4}) is similar, we use the identity
$$
G+v\frac{\partial}{\partial s}F=-v^2\frac{\partial}{\partial v}F-
-uv \frac{\partial}{\partial u}F
$$
and repeat the same steps. 

\sm

{\bf \punct Table of formulas.%
	\label{ss:list}}
 First, we present formulas for the action of the Lie algebra $\frg\frl_4$
 corresponding to $U_\sigma$, see (\ref{eq:U-sigma}). 
 Denote generators of $\frg\frl_4$ by $e_{kl}$,
 where $1\le k,l\le 4$. Denote by $\partial_{pq}$ 
 the partial derivatives $\frac{\partial}{\partial x_{pq}}$,
 where $p,q=1,2$.
  The generators $e_{kl}$ naturally split into 4 groups corresponding
  to blocks $A$, $B$, $C$, $D$ in (\ref{eq:U-sigma}). 
  
  \sm
  
  a) Generators corresponding to the block $A$ form a Lie algebra $\frg\frl_2$:
 \begin{align*}
  e_{11}= -x_{11} \partial_{11} - x_{12}\partial_{12},\qquad
  e_{12}= -x_{21} \partial_{11} - x_{22}\partial_{12},\\
  e_{21}= -x_{11} \partial_{21} - x_{12}\partial_{22},\qquad
  e_{22}= -x_{21} \partial_{21} - x_{22} \partial_{22}.
 \end{align*}
 
b) Generators corresponding to the block $D$ also form a Lie algebra $\frg\frl_2$:
\begin{align*}
e_{33}= x_{11} \partial_{11} + x_{21} \partial_{21},\qquad
e_{34}= x_{11} \partial_{12} + x_{21} \partial_{22},\\
e_{43}= x_{12} \partial_{11} + x_{22} \partial_{21},\qquad
e_{44}= x_{12} \partial_{12} + x_{22} \partial_{22}.
\end{align*}

c) Elements corresponding to the block $B$ form a 4-dimensional Abelian Lie algebra:
\begin{align*}
e_{13}= \partial_{11} + (-1 + \sigma)\frac{x_{22}}{\det X},\qquad
e_{14}= \partial_{12} - (-1 + \sigma)\frac{x_{21}}{\det X},\\
e_{23}= \partial_{21} - (-1 + \sigma)\frac{x_{12}}{\det X},\qquad
e_{24}= \partial_{22} + (-1 + \sigma)\frac{x_{11}}{\det X}.
\end{align*}

c) Elements corresponding to the block $C$ also form a 4-dimensional Abelian
Lie algebra:
\begin{align*}
e_{31} = -(x_{11}^2 \partial_{11} + x_{11} x_{12} \partial_{12} + 
x_{11} x_{21} \partial_{21} + x_{12} x_{21} \partial_{22}) + (-1 + \sigma)x_{11},
\\
e_{32} = -(x_{11} x_{21}\partial_{11} +  x_{11} x_{22} \partial_{12} + 
x_{21}^2 \partial_{21} + x_{21} x_{22} \partial_{22}) + (-1 + \sigma)x_{21},
\\
e_{41} = -(x_{11}x_{12} \partial_{11} +  x_{12}^2 \partial_{12} + 
x_{11} x_{22} \partial_{21} + x_{12} x_{22} \partial_{22}) + (-1 + \sigma)x_{12},
\\
e_{42} = -(x_{12}x_{21} \partial_{11} +  x_{12} x_{22} \partial_{12} + 
x_{21} x_{22} \partial_{21} +  x_{22}^2 \partial_{22}) + (-1 + \sigma)x_{22}.
\end{align*}

 \sm
 
Denote by $E_{kl}$ the corresponding 
  operators $E_{kl}$ on kernels $K$.
 Formulas for operators
 of groups a), b) immediately follow from the definition of the Fourier transform,
 \begin{align*}
 E_{11}&= -t \frac{\partial}{\partial t}-(1/2-\mu_1),
  \qquad &E_{12}&=\frac{\partial}{\partial t},\\
  E_{21}&=-t^2 \frac{\partial}{\partial t}+ (-1+\mu_1-\mu_2)t,
  \qquad
  &E_{22}&= t \frac{\partial}{\partial t}+(1/2 + \mu_2),
 \end{align*}
 and
 \begin{align*}
 E_{33}&=-s\frac\partial{\partial s}- (1/2+\mu_1),
 \qquad
 &E_{34}&=\frac\partial{\partial s},\\
 E_{43}&=-s^2\frac\partial{\partial s}+(-1-\mu_1+\mu_2)s,
 \qquad
& E_{44}&= s\frac\partial{\partial s}+(1/2-\mu_2).
 \end{align*}
Next,
\begin{align*}
& E_{13} := \frac{1/2 - \mu_2 + \sigma}{\mu_1 - \mu_2} \,s
 \frac\partial{\partial t}\, V_2^-
 + \frac{1/2 - \mu_1 + \sigma}{\mu_1 - \mu_2}  \Bigl(\mu_1 - \mu_2+
   s \frac\partial{\partial s}\Bigr)\,
 V_1^-
 ,
 \\
&  E_{14}=-\frac{1/2+ \sigma - \mu_1}{\mu_1 - \mu_2}\, \frac\partial{\partial s}\, V_1^- 
  -\frac{1/2 + \sigma - \mu_2}{\mu_1 - \mu_2}\,\frac\partial{\partial t}\, V_2^-
 , \\
&   E_{23}:= 
   \frac{1/2 - \mu_2 + \sigma}{\mu_1 - \mu_2}\,s\,\Bigl(-(\mu_1-\mu_2)+t
   \frac\partial{\partial t}\Bigr)\, V_2^-
   + 
   \\
   &\qquad\qquad\qquad\qquad\qquad\qquad\qquad\qquad
   \frac{1/2 - \mu_1 + \sigma}{\mu_1 - \mu_2}\,t\,\Bigl(\mu_1-\mu_2+
    s  \frac\partial{\partial s}\Bigr)\, V_1^-
  ,  \\
&  E_{24}=
  -\frac{1/2+ \sigma - \mu_1 }{\mu_1 - \mu_2} t\, \frac\partial{\partial s}\, V_1^-
   + \frac{1/2 + \mu_2 - \sigma}{\mu_1 - \mu_2}\,
  \Bigl(\mu_1-\mu_2+ t \frac\partial{\partial s}\Bigr)\,  V_2^-   
  .
\end{align*}
%
 %
%
%
and  
\begin{align*}
&E_{31} := 
 \frac{1/2 + \mu_1 + \sigma}{\mu_1 - \mu_2}\Bigl(\mu_1-\mu_2-t \frac\partial{\partial t}\Bigr)
 \,
V_1^+
- 
\frac{1/2 + \mu_2 + \sigma}{\mu_1 - \mu_2}\, t \frac\partial{\partial s} \, V_2^+
,\\
& E_{32}=\frac{1/2 +  \mu_1 +  \sigma}{\mu_1 - \mu_2}
 \frac\partial{\partial t}\, V_1^+  +\frac{1/2 + \mu_2 +  \sigma}{\mu_1 - \mu_2}  \frac\partial{\partial s}\,
 V_2^+
 ,\\
 & E_{41} = 
   \frac{1/2 + \mu_1 + \sigma}{\mu_1 - \mu_2}\,
 s \Bigl(\mu_1-\mu_2- t
  \frac\partial{\partial t}\Bigr) V_1^+
  -
 \frac{1/2 + \mu_2 + \sigma}{\mu_1 - \mu_2}\,
 t\, \Bigl(\mu_1-\mu_2+  s
  \frac\partial{\partial s}\Bigr) V_2^+ 
,  \\
& E_{42} = \frac{1/2 + \mu_1 + \sigma}{\mu_1 - \mu_2}
s \frac\partial{\partial t} V_1^+ 
 +\frac{1/2 + \mu_2 + \sigma}{\mu_1 - \mu_2} \Bigr(\mu_1-\mu_2+s\frac\partial{\partial s}\Bigr) V_2^+ 
 . 
\end{align*}

 
 
 
 {\bf \punct The case of $\GL_2(\C)$.} Notice that formulas
 in Subsections \ref{ss:R}--\ref{ss:R-formulas} for $\SL_2(\R)$
 and in Subsections \ref{ss:C-Fourier}--\ref{ss:C-formulas} are 
 very similar, except 
 the Plancherel formulas.
 
 The analog of  formula (\ref{eq:K})  is
 \begin{multline}
 K(t,s|\mu_1,\mu_1';\mu_2,\mu_2')=\\=
 \iiint\limits_{\C^3} F\bigl(
 u - t v, s u - s t v - t w, v, s v + w \bigr)\,u^{-3/2+\mu_1\|-3/2+\mu_1'}
  w^{-3/2+\mu_2\|-3/2+\mu_2'}
 \times\\\times
 \,d \Re u\,d \Im u \,d \Re v\,d \Im v \,d \Re w\,d \Im w   .
 \end{multline}
 Its derivation is based on the same change of variables 
 (\ref{eq:zamena}), its real Jacobian  is $u\ov u$.
 A further calculation one-to-one follows the calculation for $\GL_2(\R)$.

 \tt
 
 \noindent
 Math. Dept., University of Vienna; \\
 Institute for Theoretical and Experimental Physics (Moscow); \\
 MechMath Dept., Moscow State University;\\
 Institute for Information Transmission Problems.\\
 URL: http://mat.univie.ac.at/$\sim$neretin/

\end{document}